\documentclass{article}

\usepackage{amsfonts}
\usepackage{amssymb}
\usepackage{amsmath,amsthm,amscd}

\theoremstyle{plain}
\newtheorem{thm}{Theorem}

\newtheorem{conj}[thm]{Conjecture}

\theoremstyle{definition}

\newtheorem{rem}[thm]{Remark}

\begin{document}
\begin{center}
{\large A SHORT NOTE ON CONTAINMENT OF CORES

\vspace{.5cm}
 Kyungyong Lee}

\vspace{.3cm} Department of Mathematics, University of Michigan, Ann
Arbor, Michigan, USA

kyungl@umich.edu
\end{center}

\begin{abstract}
We show that cores of ideals do not preserve the inclusion.
\end{abstract}

Huneke and Swanson\cite{HunekeSwanson} raised the question of
whether, given integrally closed ideals $I \subset I'$ in a ring
$R$, it is necessarily true that $\text{core}(I) \subset
\text{core}(I')$. Hyry and Smith \cite[Corollary 5.5.1]{HyrySmith2}
gave a partial answer. The purpose of this note is to show that the
answer is no in general. Although the simplest counterexample is a
principal ideal $I$ which is contained in a power $I'$ of a maximal
ideal, we are more interested in the case when both $I$ and $I'$
have the same height.

\bigskip

Let $R=k[x,y,z,w]_{(x,y,z,w)}$ with $k$ a field of characteristic
zero  and let $\mathfrak{m}$ denote the maximal ideal of $R$. Let
$I=I_2 + \mathfrak{m}^3$, where $I_2=(x^2+yw,y^2+zw,z^2+xw)$. The
computer algebra system Macaulay shows that $I_2$ is radical. Then
it is not hard to check that $I$ is integrally closed (see
\cite[Lemma 2.1]{LazarsfeldLee}). It follows from
\cite{PoliniUlrich} or \cite{HunekeTrung} that
$\text{core}(I)=J^{n+1}:I^n$, where $J= (x^2+yw,y^2+zw,z^2+xw,w^3)$
and $n$ is the least integer such that $I^{n+1}=JI^n$. In fact $J$
is a minimal reduction of $I$ and the computer algebra system
Macaulay shows $I^2=JI$.

On the other hand, if we take $I'=\mathfrak{m}^2$ so that $I \subset
I'$, then
$\text{core}(I')=\text{core}(\mathfrak{m}^2)=\mathfrak{m}^5$ by
 \cite[Theorem 1.3]{HyrySmith} or \cite[Proposition 4.2]{CorsoPoliniUlrich}. But the computer algebra system Macaulay
shows that $\text{core}(I)=J^2:I=I^2 \not\subset \mathfrak{m}^5$,
therefore $\text{core}(I) \not\subset \text{core}(I')$. This gives a
negative answer to the above question.

\bigskip

More generally, we have the following conjecture.

\begin{conj}
Let $R=k[x_1,...,x_n]_{(x_1,...,x_n)}$ with $k$ a field of
characteristic zero and let $\mathfrak{m}$ denote the maximal ideal
of $R$. Let $I=I_d + \mathfrak{m}^{d+1}$, where $I_d$ is a complete
intersection ideal of $s$ general $d$-forms $(1\leq s< n)$. Let
$b=\lfloor\frac{dn-s+1}{d+1}\rfloor$ and $a=dn-s+1-(d+1)b$. Then
$$\emph{core}(I)=\mathfrak{m}^{a}I^{b}.$$
\end{conj}

\begin{thm}
The conjecture holds true for $d=1$.
\end{thm}
\begin{proof}
Let $X=\text{Spec} R$. By \cite[Corollary 5.3.1]{HyrySmith2}, we
have
$$\text{core}(I)= \mathcal{J}(X,n\cdot I).$$
To get the simplest log resolution $\mu:X_2 \rightarrow X_1
\rightarrow X$ of $(X,I)$, we blow up $X$ at the origin, and then
blow up the resulting surface $X_1$ along the intersection of the
exceptional divisor, say $E_1$, and the proper transform of the
linear space defined by $I_1$. Let $E_2 \subset X_2$ denote the
second exceptional divisor. Then, using the notations in
\cite[Definition 9.2.3]{Rl:pag}, we get
$$
K_{X_2/X}=(n-1)E_1 + (n+s-1)E_2, \text{ }\text{ }\text{ }\text{ }
I\cdot \mathcal{O}_{X_2}= \mathcal{O}_{X_2} (-E_1 -2E_2),
$$
so
$$\aligned\mathcal{J}(X,n\cdot I)&=\mu_* \mathcal{O}_{X_2}((n-1)E_1 + (n+s-1)E_2-n(E_1 +2E_2))\\
&=\mu_* \mathcal{O}_{X_2}(-E_1 - (n-s+1)E_2)\\
&=\mathfrak{m}^{n-s+1-2\lfloor\frac{n-s+1}{2}\rfloor}I^{\lfloor\frac{n-s+1}{2}\rfloor}.
\endaligned$$
\end{proof}

\begin{rem}
In the situation of the above conjecture, if $d\geq 2$ then
$core(I)$ and $\mathcal{J}(n\cdot I)$ do not agree in general. But
the above example leads to the following question.

If we define the \emph{deeper-core} of $I\subset R$ by
$$
\text{deeper-core}(I)=\bigcap_{J\supset I} \text{core}(J),
$$
then is
$$\text{deeper-core}(I)= \mathcal{J}(\text{Spec} R,n\cdot I)$$
always true?
\end{rem}

The author is grateful to Craig Huneke, Rob Lazarsfeld, Karen Smith
and Irena Swanson for valuable discussions and correspondence. He
would like to thank the referee for the helpful suggestion.


\begin{thebibliography}{99}
\bibitem{CorsoPoliniUlrich} Corso, A., Polini, C., Ulrich, B. (2002), Core and residual
intersections of ideals. \emph{Trans. Amer. Math. Soc.}
354:2579--2594.


\bibitem{HunekeSwanson} Huneke, C., Swanson, I. (1995), Cores of ideals in
$2$-dimensional regular local rings, \emph{Michigan Math. J.}
42:193--208.

\bibitem{HunekeTrung} Huneke, C., Trung, N. V.(2005), On the core of ideals. \emph{Compos.
Math.} 141:1--18.

\bibitem{HyrySmith2} Hyry, E., Smith, K. E. (2003), On a Non-vanishing Conjecture of
Kawamata and the Core of an Ideal, \emph{Amer. J. Math.}
125:1349--1410.


\bibitem{HyrySmith} Hyry, E., Smith, K. E. (2004), Core versus graded core, and
global section of line bundles, \emph{Trans. Amer. Math. Soc.}
356:3143--3166.

\bibitem{Rl:pag}
Lazarsfeld, R. (2004), \emph{Positivity in Algebraic Geometry II}.
volume \textbf{49} \emph{of Ergebnisse der Mathematik und ihrer
Grenzgebiete. 3. Folge. A Series of Modern Survey in Mathematics
[Results in Mathematics and Related Areas. 3rd Series. A Series of
Modern Survey in Mathematics].} Springer-Verlag, Berlin. Positivity
in vector bundles, and multiplier ideals.

\bibitem{LazarsfeldLee} Lazarsfeld, R., Lee, K. (2007), Local Syzygies of Multiplier
ideals, \emph{Invent. Math.}  167:409--418.

\bibitem{PoliniUlrich} Polini, C., Ulrich, B. (2005), A formula for the core of an
ideal. \emph{Math. Ann.} 331:487--503.



\end{thebibliography}
\end{document}